\documentclass[11pt, reqno]{amsart}
\usepackage{graphicx,amsmath,amssymb,harvard,epsfig}
\usepackage{rotating}

\usepackage{amsfonts}

\theoremstyle{plain}

\newtheorem{proposition}{Proposition}

\newcommand{\ba}{\begin{array}}
\newcommand{\ea}{\end{array}}
\newcommand{\bs}{\begin{align}\begin{split}\nonumber}
\newcommand{\bsnumber}{\begin{align}\begin{split}}
\newcommand{\es}{\end{split}\end{align}}

\newcommand{\X}{\mathcal{X}}

\numberwithin{equation}{section}

\begin{document}
\title[ ]{Improving Point and Interval Estimates of Monotone Functions
by Rearrangement}
\author[ ]{Victor
Chernozhukov$^\dag$  \ \ Iv\'an Fern\'andez-Val$^\S$ \ \ Alfred
Galichon$^\ddag$   }

\thanks{\tiny\noindent $^\dag$ Massachusetts Institute of Technology,
Department of Economics \& Operations Research Center, and
University College London, CEMMAP. E-mail: vchern@mit.edu. Research
support from the Castle Krob Chair, National Science Foundation, the
Sloan Foundation, and CEMMAP is gratefully acknowledged.}
\thanks{\tiny $\S$ Boston
University, Department of Economics. E-mail: ivanf@bu.edu. Research
support from the National Science Foundation is gratefully
acknowledged.}
\thanks{\tiny $\ddag$ Ecole Polytechnique, D$\acute{e}$partement
d'Economie. E-mail: alfred.galichon@polytechnique.edu.}
 \maketitle

\begin{abstract}
Suppose that a target function is monotonic, namely, weakly
increasing, and an  available original estimate of this target
function is not weakly increasing.  Rearrangements, univariate and
multivariate, transform the original estimate to a monotonic
estimate that always lies closer in common metrics to the target
function. Furthermore, suppose an original simultaneous confidence
interval, which covers the target function with probability at least
$1-\alpha$, is defined by an upper and lower end-point functions
that are not weakly increasing.  Then the rearranged confidence
interval, defined by the rearranged upper and lower end-point
functions, is shorter in length in common norms than the original
interval and also covers the target function with probability at
least $1-\alpha$. We demonstrate the utility of the improved point
and interval estimates with an age-height growth chart example.
\vspace{.1in}



\end{abstract}


\section{Introduction}

A common problem in statistics is the estimation of an unknown
monotonic function. Examples of monotonic functions include
biometric age-height charts, econometric demand functions, and
quantile and distribution functions. If an original, potentially
non-monotonic, estimate is available, then the rearrangement
operation from variational analysis
\cite{HLP52,lorentz:1953,villani} can be used to monotonize the
original estimate.   The rearrangement has been shown to be useful
in producing monotonized estimates of density functions
\cite{fougeres1997}, conditional mean functions
\cite{DZ2005,dette1,dette3}, and various conditional quantile and
distribution functions, see, e.g., \citeasnoun{CFG-RE-ET} and the
MIT working paper ``Quantile and Probability Curves without
Crossing'' by the authors.

In this paper, we use Lorentz inequalities and their appropriate
generalizations to show that the rearrangement of the original
estimate is not only useful for producing monotonicity, but also
always improves upon the original estimate, whenever the latter is
not monotonic. Thus, the rearranged curves are always closer to the
target curve being estimated. Furthermore, this improvement property
does not depend on the nature of the original estimate and applies
to both univariate and multivariate cases. The improvement property
of the rearrangement also extends to the construction of confidence
bands for monotone functions. We show that we can increase the
coverage probabilities and reduce the lengths of the confidence
bands for monotone functions by rearranging their upper and lower
bounds.

Monotonization has a long history in the statistical literature,
mostly in relation to isotone regression. We will not provide an
extensive literature review, but reference a few other methods most
related to the rearrangement.  \citeasnoun{Mammen1991} studies
two-step estimators, including one with smoothing in the first step
and monotonization by isotone regression in the second.
\citeasnoun{MMTW2001} show that this and many related procedures can
be recast as projections with respect to a given norm. Another
approach is the one-step procedure of \citeasnoun{Ramsay1988}, which
projects on a class of monotone spline functions called I-splines.
Later in the paper we will compare and combine these procedures with
the rearrangement.

\vspace{-.2in}

\section{Improving Point Estimates of  Monotone
Functions by Rearrangement}

\subsection{Formulation of the problem}

A basic problem in many areas of statistics is the estimation of an
unknown target function $f_0: \Bbb{R}^d \to \Bbb{R}$. Suppose we
know that  $f_0$ is monotonic, namely weakly increasing, and an
original estimate $\hat f$ is available, which is not necessarily
monotonic, but is theoretically attractive and computationally
tractable otherwise. Many common estimation methods do indeed
produce such estimates. Can they always be improved with no harm?
The answer is yes: the rearrangement method transforms the original
estimate to a monotonic estimate $\hat f^*$, and this estimate is
\textit{closer} in common metrics to the true curve $f_0$ than the
original estimate $\hat f$. Furthermore, the rearrangement is
computationally tractable, and thus preserves the appeal of the
original estimates.

Estimation methods used in regression analysis can be grouped into
global methods and local methods. An example of a global method is
the series estimator of $f_0$ taking the form $\hat f(x) =
P_{k_n}(x)'\hat b,$ where $P_{k_n}(x)$ is a $k_n$-vector of suitable
transformations of the variable $x$, such as B-splines, polynomials,
and trigonometric functions, and
$$
\hat b = \arg \min_{b \in \Bbb{R}^{{k}_n}} \sum_{i=1}^n \rho\{Y_i -
P_{k_n}(X_i)' b\},
$$
where $\{(Y_i, X_i), i=1,\ldots,n\}$ denotes the data.  In
particular, using the square loss $\rho(u) = u^2$ produces estimates
of the conditional mean of $Y_i$ given $X_i$
\cite{gallant:fourier,andrews:series,Stone:series,newey:series},
while using the asymmetric absolute deviation loss $\rho(u) = \{u -
1(u<0)\} u$ produces estimates of the conditional $u$-quantile of
$Y_i$ given $X_i$ \cite{koenker:1978,portnoy:splines,he:series}. The
series estimates $x \mapsto \hat f(x)= P_{k_n}(x)'\hat b$ are widely
used in data analysis due to their desirable approximation and
theoretical properties, and computational tractability. However,
they need not be monotone, unless explicit constraints are added
\cite{matzkin:handbook,silvapulle:book,koenker:inequality}.

Examples of local methods include kernel and local polynomial
estimators. A kernel estimator takes the form
$$
\hat f(x) = \arg \min_{b \in {\Bbb{R}}} \ \sum_{i=1}^n w_i
 \rho(Y_i - b), \ \  w_i  = K\left(\frac{X_i
 -x}{h}\right),
$$
where the loss function $\rho$ plays the same role as above, $K(u)$
is a multivariate kernel function, and $h> 0$ is a vector of
bandwidths \cite{wand:jones,silverman:book}. The resulting estimate
$x \mapsto \hat f(x)$ need not be monotone. \citeasnoun{dette1} show
that the rearrangement transforms the kernel estimate into a
monotonic one. We further show here that the rearranged estimate
necessarily improves upon the original estimate, whenever the latter
is not monotonic. Local polynomial regression is a related local
method \cite{chaudhuri:1991,fan:book}. In particular, the local
linear estimator takes the form
$$
\{\hat f(x), \hat d(x)\} = \underset{b \in {\Bbb{R}}, c \in
\Bbb{R}^d}{\text{argmin}} \ \sum_{i=1}^n w_i \rho\{Y_i - b - c'(X_i
-x)\}^2, \ \  w_i  = K\left(\frac{X_i -x}{h}\right).
$$
The resulting estimate $x \mapsto \hat f(x)$, while theoretically
attractive and computationally tractable, may also be non-monotonic,
as illustrated in Section 4.

\subsection{The rearrangement and its estimation property: the univariate case}
In what follows, let $\mathcal{X}$ be a compact interval; without
loss of generality we take $\mathcal{X}=[0,1]$. Let $f$ be a
measurable function mapping $\mathcal{X}$ to $K$, a bounded subset
of $\Bbb{R}$.  The increasing rearrangement  $f^*$ of $f$ is  the
quantile function of the random variable $f(X)$ when $X\sim U(0,1)$,
that is,
$$f^*(x) = \inf \left\{ y \in \Bbb{R}: \int_{\mathcal{X}}
1\{ f(u) \leq y \} du  \geq x \right\}.
$$
The rearrangement operator simply transforms a function $f$ to its
quantile function $f^*$. For computing purposes when $f$ is
continuous, we can think of the rearrangement as a sorting
operation: given values of the function $f$ evaluated at $x$ in a
fine enough net of equidistant points, we simply sort the values in
increasing order to create the sorted, i.e., rearranged, function.

\begin{proposition}\label{betterapproximation}
Let the target $f_0: \mathcal{X} \to K$ be a weakly increasing measurable
function in $x$, and $\hat f:
\mathcal{X} \to K$ be another measurable function, an initial
estimate of  $f_0$.

\noindent 1. For any $p \in [1,\infty]$, the rearrangement of $\hat
f$, denoted $\hat f^*$, weakly reduces the estimation error:
 \begin{equation}\label{P11}
\left\{\int_{\mathcal{X}} \left| \hat f^*(x) - f_0(x) \right|^p d x
\right\}^{1/p} \leq \left\{ \int_{\mathcal{X}} \left| \hat f(x) -
f_0(x) \right|^p d x \right\}^{1/p}.
 \end{equation}

\noindent 2. Suppose that there exist regions $\X_0$ and $\X'_0$,
each of measure greater than $\delta>0$, such that for all $x \in
\X_0$ and $x' \in \X'_0$ we have that (i) $x'> x$, (ii) $\hat f(x) >
\hat f(x') + \epsilon$, and (iii) $f_0(x')
> f_0(x) + \epsilon$, for some $\epsilon >0$.  Then the
gain in the quality of estimation is strict for $p \in
(1,\infty)$. Namely, for any $p \in (1,\infty)$,
 \begin{equation}\label{P12}
\left\{\int_{\mathcal{X}} \left| \hat f^*(x) -f_0(x) \right|^p d x
\right\}^{1/p} \leq \left\{ \int_{\mathcal{X}} \left| \hat f(x)
-f_0(x)\right|^p d x  - \delta \eta_p \right\}^{1/p},
\end{equation}
where $\eta_p = \inf\{ |v - t'|^p + |v' - t|^p - |v-t|^p - |v'-t'|^p
\}>0$, with the infimum taken over all $v, v', t, t'$ in the set $K$
such that $v' \geq v + \epsilon$ and $t' \geq t + \epsilon$.

\end{proposition}

Proposition 1 establishes that the rearranged estimate $\hat f^*$
has a smaller, often strictly smaller, estimation error in the $L_p$
norm than the original estimate whenever the latter is not monotone.
This very useful and generally applicable property is independent of
the sample size and of the way the original estimate $\hat f$ is
obtained. As follows from (\ref{P12}), the reduction in estimation
error is strict for $L^p$ norms with $p \in (1, \infty)$ if the
original estimate $\hat f$ is decreasing on a subset of
$\mathcal{X}$ having positive measure, while the target function
$f_0$ is increasing on this subset. If $f_0$ is constant, then there
is no reduction in estimation error; that is, the inequality
(\ref{P11}) becomes an equality, since the random variables $\hat
f^*(X)$ and  $\hat f(X)$ share the same quantile function $\hat f^*$
and hence the same distribution function, and $f_{0}(X)$ is
constant.

The weak inequality (\ref{P11}) is a direct, yet important,
consequence of the classical rearrangement inequality due to
\citeasnoun{lorentz:1953}: let $q$ and $g$ be two functions mapping
$\mathcal{X}$ to $K$, and $q^{*}$ and $g^{*}$ be their corresponding
increasing rearrangements, then $ \int_{\mathcal{X}} L\{q^*(x),
g^*(x)\} d x  \leq \int_{\mathcal{X}} L\{q(x), g(x) \} dx, $ for any
submodular discrepancy function $L: \Bbb{R}^2 \mapsto \Bbb{R}_+$. We
set $q = \hat f$, $q^* = \hat f^*$, $g = f_0$, and $g^* = f_0^*$. In
our case $ f_0^* = f_0$ almost everywhere, that is, the target
function is its own rearrangement. Further, recall that $L$ is
submodular if for each pair of vectors $(v,t)$ and $(v',t')$ in
$\Bbb{R}^2$, we have that
\begin{equation}\label{submodular}
L(v\wedge v', t \wedge t') + L(v \vee v', t \vee t') \leq L( v, t) +
L(v', t').
\end{equation}
In other words, a function $L$ measuring the discrepancy between
pairs of vectors is submodular if co-monotonization of the pair
reduces the discrepancy. When the function $L$ is smooth,
submodularity is equivalent to   $\partial^2 L(v,t)/ (\partial v
\partial t)  \leq 0$ holding for each $(v, t)$ in $\Bbb{R}^2$. Thus,
for example, power functions $L(v,t) = |v-t|^p$ for $p \in
[1,\infty)$ and many other loss functions are submodular.  The weak
inequality (\ref{P11}) then follows.

\subsection{The rearrangement and its estimation
property: the multivariate case}

In this section we consider multivariate functions $f: \X^d \to K$,
where $\X^d = [0,1]^d$ and $K$ is a bounded subset of $\Bbb{R}$. The
notion of monotonicity we seek to impose on $f$ is the following: we
say that the function $f$ is weakly increasing in the vector $x$
if $f(x') \leq f(x)$ whenever $x' \leq x$ (componentwise). 
In what follows, we use $f(x_j, x_{-j})$ to denote the dependence of
$f$ on $x_j$, and all other arguments, $x_{-j}$, that exclude $x_j$.
The notion of monotonicity above is equivalent to the requirement
that for each $j$ in $1,\ldots,d$ the mapping $x_j \mapsto f(x_j,
x_{-j})$ is weakly increasing in $x_j$, for each $x_{-j}$ in
$\X^{d-1}$.

Define the rearrangement operator $R_j$ and the rearranged function
$f_j^*$ with respect to $x_j$ as
$$ f_j^*(x) = R_j f(x) = \inf \left \{ y: \left [\int_{\mathcal{X}}
1 \{  f(x_j', x_{-j}) \leq y \} d x_j' \right] \geq x_j \right \}.$$
This is the one-dimensional increasing rearrangement applied to the
one-dimensional function $x_j \mapsto f(x_j, x_{-j})$, holding the
other arguments $x_{-j}$ fixed. The rearrangement is applied for
every value of the other arguments $x_{-j}$.

Let $\pi = (\pi_1,\ldots,\pi_d)$ be an ordering, i.e., a permutation,
of the integers $1,\ldots,d$.  Let us define the $\pi$-rearrangement
operator $R_{\pi}$ and the $\pi$-rearranged function $f_\pi^*$ as
$
f_{\pi}^* = R_\pi f = R_{\pi_1} \ldots  R_{\pi_d}
 f.$
For any ordering $\pi$, the $\pi$-rearrangement operator rearranges
the function with respect to all of its arguments.  As shown below,
the resulting function $f_\pi$ is weakly increasing in $x$. In
general, two different orderings $\pi$ and $\pi'$ of $1,\ldots,d$
can yield different rearranged functions $f_{\pi}^*$ and
$f_{\pi'}^*$.  To resolve the conflict among rearrangements done
with different orderings, we may consider averaging among them:
letting $\Pi$ be any finite collection of orderings $\pi$, we can
define the average rearrangement as
$$
f^*  = \frac{1}{|\Pi|} \sum_{\pi \in \Pi} f_\pi^*,
$$
where $|\Pi|$ denotes the number of elements in the set of orderings
$\Pi$. \citeasnoun{dette3} also proposed averaging all the possible
orderings of a related smoothed procedure in the context of monotone
conditional mean estimation. As shown below, the estimation error of
the average rearrangement is weakly smaller than the average of
estimation errors of individual $\pi$-rearrangements.

The following proposition describes the properties of multivariate
$\pi$-rearrangements:

\begin{proposition}\label{betterapproximationRn}  Let the target function
$f_0: \X^d \to K$ be weakly increasing and measurable in $x$. Let
$\hat f: \X^d \to K$ be a measurable function that is an initial
estimate of  $f_0$. Let $\bar f: \X^d \to K$ be
another estimate of $f_0$, which is measurable in $x$, including,
for example, a rearranged $\hat f$ with respect to some of the
arguments.

 1.  For each ordering $\pi$ of $1,\ldots,d$, the
$\pi$-rearranged estimate $\hat f^*_\pi$ is weakly increasing.   Moreover, $\hat f^*$, an average of $\pi$-rearranged
estimates, is weakly increasing.

 2. (a) For any $j$ in $1,\ldots, d$ and any $p$
in $[1,\infty]$, the rearrangement of $\bar f$ with respect to the
$j$-th argument produces a weak reduction in the estimation
error:
 \begin{align*}\begin{split}
\left\{\int_{\X^d} | \bar f_j^*(x) - f_0(x) |^p dx\right \}^{1/p}
\leq \left\{\int_{\X^d} | \bar f(x) - f_0(x) |^p dx\right \}^{1/p}.
\end{split}\end{align*}

(b) A $\pi$-rearranged estimate $\hat f^*_\pi$ of $\hat f$ weakly
reduces the estimation error of $\hat f$:
 \bsnumber\label{P22b}
\left\{\int_{\X^d} | \hat f^*_{\pi}(x) - f_0(x) |^p dx\right\}^{1/p}
\leq \left\{\int_{\X^d} | \hat f(x) - f_0(x) |^p dx\right\}^{1/p}.
 \end{split}\end{align}

3. Suppose that there exist subsets $\X_j \subset \X$ and $\X_j' \subset
\X$, each of measure greater than $\delta>0$, and a subset $\X_{-j} \subseteq
\X^{d-1}$, of measure $\nu>0$, such that for all  $x= (x_j, x_{-j})$
and $x' = (x'_j, x_{-j})$, with $x_j' \in \X_j'$, $x_j \in \X_j$,
$x_{-j} \in \X_{-j}$, we have that (i) $x_j' > x_j$, (ii) $\bar f(x)
> \bar f(x') + \epsilon$, and (iii) $f_0(x') > f_0(x) + \epsilon$,
for some $\epsilon>0$.

(a) Then, for any $p \in (1, \infty)$,
 \begin{align*}\begin{split}
\left\{\int_{\X^d} | \bar f_j^*(x) - f_0(x) |^p dx\right\}^{1/p}
\leq \left\{\int_{\X^d} | \bar f(x) - f_0(x) |^p dx - \eta_p \delta
\nu \right\}^{1/p},
\end{split}\end{align*}
where $\eta_p = \inf\{ |v - t'|^p + |v' - t|^p - |v-t|^p - |v'-t'|^p
\}>0$, with the infimum taken over all $v, v', t, t'$ in the set $K$
such that $v' \geq v + \epsilon$ and $t' \geq t + \epsilon$.

(b) Further, for an ordering $\pi  = (\pi_1,\ldots,
\pi_k,\ldots,\pi_d)$ with $\pi_k = j$, let $\bar f$ be a partially
rearranged function, $\bar f = R_{\pi_{k+1}} \ldots R_{\pi_d} \hat
f$ (for $k =d$ we set $\bar f = \hat f$). If the function $\bar f$
and the target function $f_0$ satisfy the condition stated above,
then, for any $p \in (1,\infty)$,
 \bsnumber\label{P23b}
\left\{\int_{\X^d}  | \hat f_\pi^*(x) - f_0(x) |^p dx\right\}^{1/p}
\leq \left\{\int_{\X^d} | \hat f(x) - f_0(x) |^p dx - \eta_p \delta
\nu \right\}^{1/p}.
\end{split}\end{align}

 4.  The estimation error of an average
rearrangement is weakly smaller than the average estimation error of
the individual $\pi$- rearrangements: for any $p \in [1, \infty]$,
 \begin{align*}\begin{split}
\left\{\int_{\X^d}  | \hat f^*(x) - f_0(x) |^p dx\right\}^{1/p}
\leq \frac{1}{|\Pi|} \sum_{\pi \in \Pi} \left\{\int_{\X^d} | \hat
f_\pi^*(x) - f_0(x) |^p dx\right\}^{1/p}.
\end{split}\end{align*}
\end{proposition}

Proposition 2 generalizes Proposition 1 to the multivariate case,
also demonstrating several features unique to the multivariate case.
We see that the $\pi$-rearranged functions are monotonic in all of
the arguments. \citeasnoun{dette3}, using a different argument,
showed that their related smoothed procedure for conditional mean
functions is monotonic in both arguments for the bivariate case in
large samples. The rearrangement along any argument improves the
estimation properties. Moreover, the improvement is strict when the
rearrangement with respect to a $j$-th argument is performed on an
estimate that is decreasing in the $j$-th argument, while the target
function is increasing in the same $j$-th argument, in the sense
precisely defined in the proposition. Averaging different
$\pi$-rearrangements is better on average than using a single
$\pi$-rearrangement chosen at random.

\subsection{Discussion}
%
\begin{figure}
\includegraphics[width = \textwidth, height = .5\textwidth]{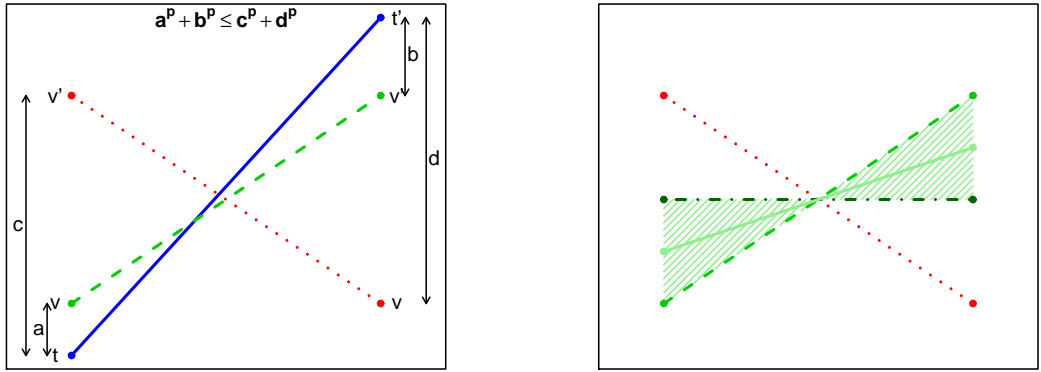}
\caption{Geometric illustration for the proof of Proposition 1 (left
panel) and comparison to isotonic regression (right panel).  The
solid dark line is the target function $f_0$, the dotted line is the
original estimate $\hat{f}$, the dashed line is the rearranged
estimate $\hat{f}^{\ast}$, the dotted-dashed line is the isotonized
estimate $\hat{f}^{I}$, and the solid light line is the average of
the rearranged and isotonized estimates $\hat{f}^{1/2}$. In the left
panel $L(v,t) = a^p$, $L(v',t) = c^p$, $L(v', t') = b^p$, and $L(v,
t') = d^p$.} \label{Fig:prop1}
\end{figure}

Here we informally explain why rearrangement provides the
improvement property and compare rearrangement to isotonization.

We begin by noting that the proof of the improvement property
can be first reduced to the case of step functions or,
equivalently, functions with a finite domain, and then to the case
of functions with a two-point domain. The improvement
property for such functions then follows from the
submodularity property (\ref{submodular}). In the left panel of
Figure \ref{Fig:prop1} we illustrate  this geometrically by
plotting the original estimate $\hat f$, the rearranged estimate
$\hat f^*$, and the true function $f_0$.  In this example, the
original estimate is decreasing and hence violates the monotonicity
requirement. We see that the two-point rearrangement co-monotonizes
$\hat f^*$ with $f_0$ and thus brings $\hat f^*$ closer to $f_0$.
Also, we can view the rearrangement as a projection on the set of
weakly increasing functions that have the same distribution as the
original estimate $\hat f$.

In the right panel of Fig. \ref{Fig:prop1} we plot both the
rearranged and isotonized estimates. The isotonized estimate $\hat
f^I$ is a projection of the original estimate $\hat f$ on the set of
weakly increasing functions, that only preserves the mean of the
original estimate.  We can compute the two values of the isotonized
estimate $\hat f^I$ by assigning to both the average of the two
values of the original estimate $\hat f$, whenever the latter
violate the monotonicity requirement, and leaving the original
values unchanged otherwise.  In our example in Fig. \ref{Fig:prop1}
this produces a flat function $\hat f^I$. This pool adjacent
violators procedure extends to domains with more than two points by
applying the procedure iteratively to any pair of points at which
 monotonicity is violated \cite{PAVA1955}.

Using the computational definition of isotonization, one can show
that, like rearrangement, isotonization also improves upon the
original estimate, for any $p \in[1,\infty]$:
\begin{equation*}
\left\{\int_{\mathcal{X}} | \hat f^{I}(x) - f_0(x) |^p d x  \right\}^{1/p}   \leq  \left\{\int_{\mathcal{X}} | \hat f(x) - f_0(x) |^p d
x  \right\}^{1/p},
\end{equation*}
see, e.g., \citeasnoun{barlow:book}. Therefore, it follows that any
function $\hat f^{\lambda}$ in the convex hull of the rearranged and
isotonized estimate both (1) monotonizes and (2) improves upon the
original estimate $\hat f$, that is, for any $p \in[1,\infty]$ and
$\lambda \in [0,1]$,
\begin{equation*} \left\{\int_{\mathcal{X}}
| \hat f^{\lambda}(x) - f_0(x) |^p d x \right\}^{1/p}  \leq    \ \
\left\{ \int_{\mathcal{X}} | \hat f(x) - f_0(x) |^p d x
\right\}^{1/p},
\end{equation*}where $\hat f^{\lambda} = \lambda \hat f^* + (1 - \lambda)
\hat f^I$.  The first property is obvious and the second follows
from homogeneity and subadditivity of norms. By induction on the
dimension, the improvement property extends to the sequential
multivariate isotonization and to its convex hull with the
sequential multivariate rearrangement.

Thus, we see that a rather rich class of procedures both monotonizes
the original estimate and reduces the distance to the true target
function. However, there is no single best distance-reducing
monotonizing procedure. Indeed, whether the rearranged estimate
$\hat f^*$ approximates the target function better than the
isotonized estimate $\hat f^I$ depends on how steep or flat the
target function is. We illustrate this point using the example
plotted in the right panel of Fig. \ref{Fig:prop1}:  consider any
increasing target function taking values in the shaded area between
$\hat f^*$ and $\hat f^I$, and also the function $\hat f^{1/2}$, the
average of the isotonized and the rearranged estimate, that passes
through the middle of the shaded area. Suppose first that the target
function is steeper than $\hat f^{1/2}$, then $\hat f^*$ has a
smaller estimation error than $\hat f^I$. Now suppose instead that
the target function is flatter than $\hat f^{1/2}$,  then $\hat f^I$
has a smaller estimation error than $\hat f^*$. It is also clear
that, if the target function is neither very steep nor very flat,
$\hat f^{1/2}$ can outperform either $\hat f^*$ or $\hat f^I$. Thus,
in practice we can choose rearrangement, isotonization, or, some
combination of the two, depending on our beliefs about how steep or
flat the target function is in a particular application.

\vspace{-.2in}

\section{Improving Interval Estimates of Monotone Functions by Rearrangement}

In this section we propose to directly apply the rearrangement,
univariate and multivariate, to simultaneous confidence intervals
for monotone functions. We show that our proposal will necessarily
improve the original intervals by decreasing their length while
retaining the same or greater coverage level.

Suppose that we are given an initial simultaneous confidence
interval
\begin{equation}\label{CI1}
[\ell, u]= \{ [\ell(x), u(x)], x \in \mathcal{X}^d \},
\end{equation}
where $\ell$ and $u$ are the lower and upper end-point
functions such that $\ell \leq u$ on $\mathcal{X}^d$, that is, $\ell(x) \leq u(x)$ for all $x \in \mathcal{X}^d$. We further suppose that the confidence interval $[\ell, u]$ has either the exact or the asymptotic \emph{confidence property} for the estimand function $f$, namely, for a given $\alpha \in (0,1)$,
\begin{equation}\label{CP}
\text{pr}_P \{f \in [\ell, u] \} \geq 1- \alpha,
\end{equation}
for all probability measures $P$ in some set $\mathcal{P}_n$
containing the true probability measure $P_0$. The statement $f \in
[\ell, u]$ means that  $\ell(x) \leq f(x) \leq u(x)$ for all $x \in
\mathcal{X}^d$. We assume that property (\ref{CP}) holds either in
the finite sample sense, that is, for the given sample size $n$, or
in the asymptotic sense, that is, for all but finitely many sample
sizes $n$ \cite{LR2005}.

A common confidence interval for functions specifies
\begin{equation}
\label{CI2} \ell(x) = \hat f(x) - \hat c s(x),  \ \ \ \ \ u(x) =
\hat f(x) + \hat c s(x),
\end{equation}
where $\hat f(x)$ is a point estimate, $s(x)$ is the standard error
of the point estimate, and $\hat c$ is a critical value chosen to
attain the confidence property (\ref{CP}).  \citeasnoun{Wasserman06}
provides an excellent overview of methods for constructing the
critical value. The problem with such confidence intervals, as with
the point estimates themselves, is that they need not be monotonic.
Indeed, the end-point functions (\ref{CI2}) need not be monotonic,
so the confidence interval may contain non-monotone functions
excludable from it. Accordingly we can intersect the interval with
the set of monotone functions to reduce its length without affecting
its coverage level. In some cases, however, the initial interval may
not contain any monotone function and the resulting intersected
interval is empty, due, for example, to misspecification.

We say that confidence intervals are misspecified or incorrectly
centered if the estimand $f$, being covered by $[\ell, u]$ in
(\ref{CP}), is not equal to the weakly increasing target function
$f_0$, so that $f$ may not be monotone.  Incorrect centering is
rather common both in parametric and non-parametric estimation. In
parametric estimation correct centering of confidence intervals
requires perfect specification of functional forms,  whereas in
nonparametric estimation correct centering requires the so-called
undersmoothing; both are difficult.  In real applications with many
regressors, researchers tend to use oversmoothing rather than
undersmoothing. In a recent development, \citeasnoun{GW08} provide
some formal justification for oversmoothing: targeting inference on
functions $f$, that represent various smoothed versions of $f_0$ and
thus summarize features of $f_0$, may be desirable to make inference
more robust, or, equivalently, to enlarge the class of
data-generating processes $\mathcal{P}_n$ for which (\ref{CP})
holds. Regardless of the reasons for why the confidence intervals
may target $f$ instead of $f_0$, our procedures will work for
inference on the monotonized, hence improved, version $f^*$ of $f$.

Our proposal for improved interval estimates is to rearrange the
entire simultaneous confidence interval into a monotonic interval
\begin{equation}\label{RCI}
[\ell^*, u^*] = \{ [\ell^*(x), u^*(x)], x \in \mathcal{X}^d  \},
 \end{equation}
where the lower and upper end-point functions $\ell^*$ and $u^*$ are
the increasing rearrangements of the original end-point functions
$\ell$ and $u$. In the multivariate case, we use the symbols
$\ell^*$ and $u^*$ to denote either multivariate
$\pi$-rearrangements $\ell_\pi^*$ and $u_\pi^*$ or average
multivariate rearrangements $\ell^*$ and $u^*$, whenever we do not
need to emphasize specifically the dependence on $\pi$.

 The following proposition describes the properties of the rearranged confidence intervals.

\begin{proposition} Let $[\ell, u]$ in (\ref{CI1}) be the original confidence interval satisfying the property
(\ref{CP}) for the estimand function $f : \mathcal{X}^d \mapsto K$
and let the rearranged confidence interval $[\ell^*, u^*]$ be
defined as in (\ref{RCI}).

1.  The interval $[\ell^*, u^*]$ is weakly
increasing and non-empty, in the sense that the end-point functions
$\ell^*$ and $u^*$ are weakly increasing on $\mathcal{X}^d$ and
satisfy $\ell^* \leq u^*$ on $\mathcal{X}^d$.  Moreover, the event
that  $f \in [\ell, u]$ implies the event that $f^* \in [\ell^*, u^*].$
In particular, under the correct specification, when $f$ equals a
weakly increasing target function $f_0$, we have that $f=f^*=f_0$,
so that $f_0 \in [\ell, u]$ implies $f_0 \in [\ell^*, u^*].$
Therefore, $[\ell^*, u^*]$ covers $f^*$, which is equal to $f_0$
under the correct specification,  with a probability that is greater
or equal to the probability that  $[\ell, u]$ covers $f$.

2.  The interval $[\ell^*, u^*]$ is weakly
shorter than $[\ell, u]$ in the
 $L^p$ length: for each $p \in [1,\infty]$,
\begin{equation}\label{shorter1}
\left\{\int_{\mathcal{X}^d}  \Big | \ell^*(x) - u^*(x) \Big | ^p
dx\right\}^{1/p} \leq \left\{\int_{\mathcal{X}^d} \Big | \ell(x) -
u(x) \Big |^p dx\right\}^{1/p}.
\end{equation}

3.  In the univariate case, suppose that there exist subsets $\X_0
\subset \X$ and $\X'_0 \subset \X$, each of measure greater than
$\delta>0$ such that for all $x' \in \X'_0$ and $x \in \X_0$,  we
have that $x' > x$, and either (i) $\ell(x) > \ell (x') + \epsilon$,
and $ u (x') > u (x) + \epsilon$, for some $\epsilon>0$ or (ii)
$\ell (x') > \ell (x) + \epsilon$ and $ u (x) >  u (x') + \epsilon$,
for some $\epsilon>0$.  Then, for any $p \in (1,\infty)$,
\begin{equation*}
\left\{\int_{\mathcal{X}}  \Big |
\ell^*(x) - u^*(x) \Big | ^p dx\right\}^{1/p} \leq \left\{\int_{\mathcal{X}} \Big | \ell(x) - u(x) \Big |^p - \eta_p \delta
\right\}^{1/p},
\end{equation*}
where $\eta_p = \inf\{ |v - t'|^p + |v' - t|^p - |v-t|^p - |v'-t'|^p
\}>0$, where the infimum is taken over all $v, v', t, t'$ in $K$
such that $v' \geq v + \epsilon$ and $t' \geq t + \epsilon$.

In the multivariate case with $d\geq 2$, for an ordering $\pi  =
(\pi_1,\ldots, \pi_k,\ldots,\pi_d)$ of integers $\{1,\ldots,d\}$
with $\pi_k = j$, let $\bar g$ denote the partially rearranged
function,  $\bar g = R_{\pi_{k+1}}  \ldots R_{\pi_d}  \hat g$, where
for $k =d$ we set $\bar g = \hat g$.    Suppose there exist subsets
$\X_j \subset \X$ and $\X_j' \subset \X$, each of measure greater
than $\delta>0$, and a subset $\X_{-j} \subseteq \X^{d-1}$, of
measure $\nu>0$, such that for all  $x= (x_j, x_{-j})$ and $x' =
(x'_j, x_{-j})$, with $x_j' \in \X_j'$, $x_j \in \X_j$, $x_{-j} \in
\X_{-j}$, we have that (i) $x_j' > x_j$, and either (ii) $\bar
\ell(x) > \bar \ell (x') + \epsilon$, and $ \bar u (x') > \bar u (x)
+ \epsilon$, for some $\epsilon>0$ or (iii) $\bar \ell (x') > \bar
\ell (x) + \epsilon$ and $ \bar u (x) > \bar u (x') + \epsilon$, for
some $\epsilon>0$. Then, for any $p \in (1,\infty)$ and $\eta_p>0$
defined as above
\begin{equation*}
\left\{\int_{\mathcal{X}^d}  \Big | \ell^*_\pi(x) - u^*_\pi(x) \Big
| ^p dx\right\}^{1/p} \leq \left\{\int_{\mathcal{X}^d} \Big |
\ell(x) - u(x) \Big |^p - \eta_p \delta \nu \right\}^{1/p}.
\end{equation*}

\end{proposition}

Proposition 3 shows that the rearranged confidence intervals are
weakly shorter than the original confidence intervals, and also
qualifies when the rearranged  confidence intervals are strictly
shorter. In particular, in the univariate case the inequality
(\ref{shorter1}) is necessarily strict for $p\in (1, \infty)$ if
there is a region of positive measure in $\mathcal{X}$ over which
the end-point functions $\ell$ and $u$ are not comonotonic.  This
weak shortening result follows for univariate cases directly from
the Lorentz (1953) inequality, and the strong shortening by its
strengthening.  The shortening results for the multivariate case
follow by induction on the dimension. Moreover, the
order-preservation property of the univariate and multivariate
rearrangements, demonstrated in the proof, implies that the
rearranged confidence interval $[\ell^*, u^*]$ has a weakly higher
coverage than the original confidence interval $[\ell, u]$. We do
not quantify strict improvements in coverage, but demonstrate them
through the examples in the next section.

Our idea of directly monotonizing the interval estimates also
applies to other monotonization procedures. Indeed,  the proof of
Proposition 3 reveals that  part 1 applies to any
order-preserving monotonization operator $T$, such that
\begin{equation}\label{order-preserve}
 g \leq m  \text{ implies } T g \leq T m.
\end{equation}
Furthermore, part 2 of Proposition 3 on the weak shortening of the
confidence intervals applies to any distance-reducing operator $T$
such that
 \begin{equation}\label{distance-reduce}
\left\{ \int_{\mathcal{X}^d}| T \ell(x) - T u(x) |^p d x
\right\}^{1/p} \leq \left\{ \int_{\mathcal{X}^d}| \ell(x) - u(x) |^p
d x \right\}^{1/p}.
 \end{equation}
Rearrangements are instances of operators that have properties
(\ref{order-preserve}) and (\ref{distance-reduce}). Isotonization is
another important instance \cite{RWD88}. Moreover, convex
combinations of order-preserving and distance-reducing operators,
such as the average of rearrangement and isotonization, also have
properties (\ref{order-preserve}) and (\ref{distance-reduce}).
\vspace{-.2in}
\section{Illustrations}

\subsection{An empirical illustration with
age-height reference charts}

In this section we provide an empirical application to biometric
age-height charts. We show how the rearrangement monotonizes and
improves various nonparametric point and interval estimates for
functions.

Since their introduction by Quetelet in the 19th century, reference
growth charts have become common tools to assess an individual's
health status. These charts describe the evolution of individual
anthropometric measures, such as height, weight, and body mass
index, across different ages. See \citeasnoun{cole} for a classical
work on the subject, and \citeasnoun{koenker:charts} for a recent
analysis from a quantile regression perspective and additional
references. Here we consider an application of the rearrangement and
other related methods to the estimation of growth charts for height.
This makes sense since an individual's height should follow an
increasing relationship with age up to adulthood. Our data consist
of repeated cross sectional measurements of height in centimeters
and age in months from the 2003-2004 US National Health and
Nutrition Survey, and is further restricted to the subsample of
US-born white males aged 2-20 to avoid other confounding factors,
giving us a sample of 533 observations.

Let $Y$ and $X$ denote height and age, respectively. Let $E[Y\mid
X=x]$ denote the conditional expectation  of $Y$ given $X=x$, and
$Q_Y[u\mid X=x]$ denote the conditional $u$-th quantile of $Y$ given
$X=x$, where $u$ is the quantile index. The target functions of
interests are the conditional expectation function, $x \mapsto
E[Y\mid X=x]$, the conditional quantile functions for several
quantile indices, $x\mapsto Q_Y[u\mid X=x]$, for $u =5\%$, $50\%$,
and $95\% $, and the entire conditional quantile process  for height
given age, $(u,x)\mapsto Q_Y[u\mid X=x]$.   The monotonicity
requirements for these target functions are the following: the first
two should be increasing in age $x$, and the third should be
increasing in both age $x$ and the quantile index $u$.

We estimate the target functions using non-parametric ordinary least
squares or quantile regression and then rearrange the estimates to
satisfy the monotonicity requirements. We consider kernel, local
linear, regression splines, and Fourier series methods.  For the
kernel and local linear methods, we choose a bandwidth of one year
and a box kernel. For the regression splines method, we use cubic
B-splines with a knot sequence $\{ 3, 5, 8, 10, 11.5, 13, 14.5, 16,
18 \}$ \cite{koenker:charts}. For the Fourier method, we employ four
sines and four cosines. For the estimation of the conditional
quantile process, we use  $\{ 0.005, 0.010, \ldots, 0.995 \}$ as a
net of quantile indices.

\begin{figure}
\includegraphics[width = \textwidth, height = \textwidth]{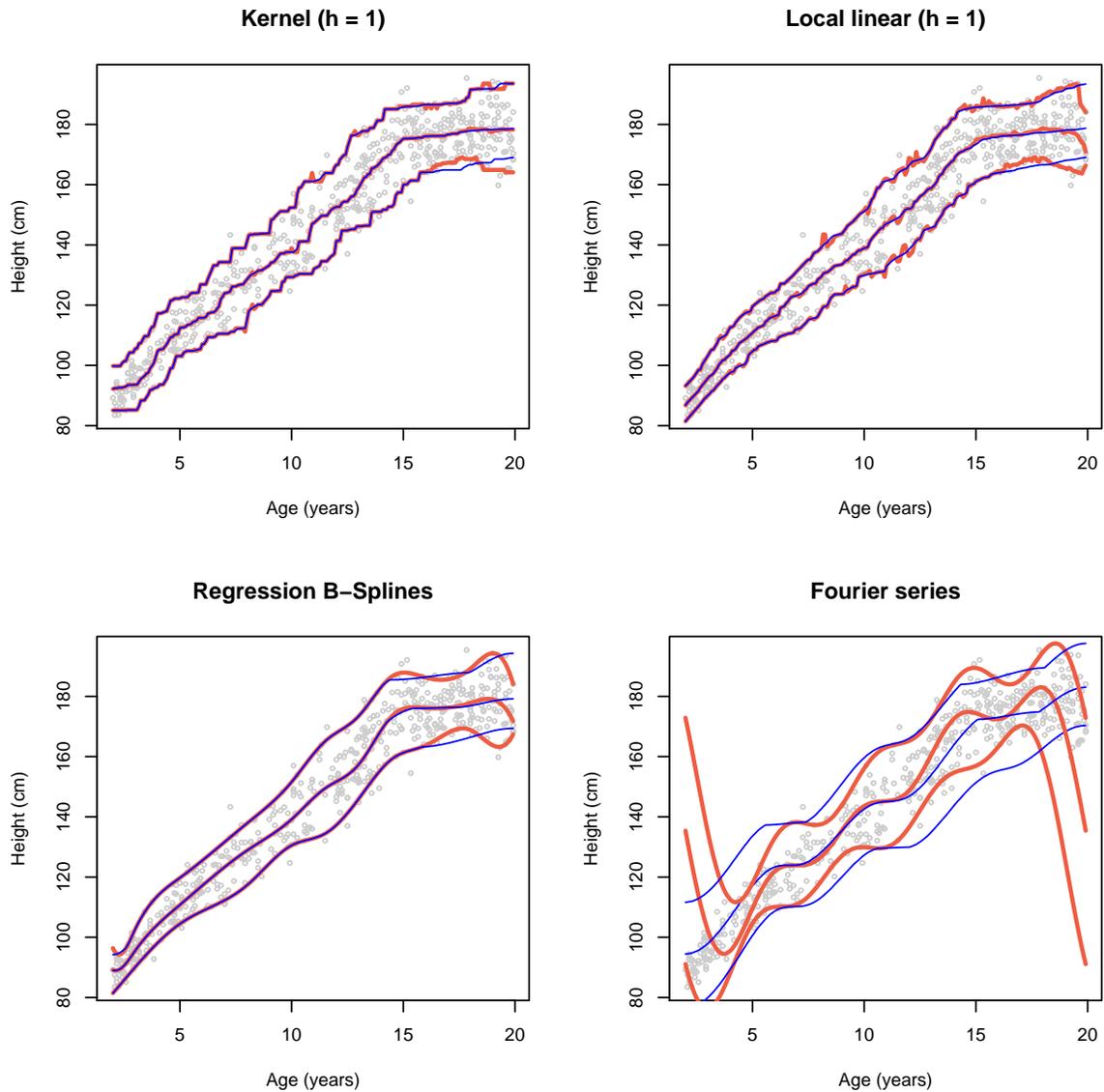}
\caption{Estimates of the 5\%, 50\%, and 95\% conditional quantile
functions of height given age and their increasing rearrangements,
obtained by kernel, local linear, cubic B-splines series,  and
Fourier series regression. Light thick lines are the original
estimates and dark thin lines are the rearranged estimates.}
\label{Fig:2}
\end{figure}

Figure \ref{Fig:2} shows the original and
rearranged estimates of the conditional quantile functions for the
different methods. All the estimated curves have trouble capturing
the slowdown in the growth of height after age fifteen and yield
non-monotonic curves for the highest values of age. The Fourier
series performs particularly poorly in approximating the aperiodic
age-height relationship and has many non-monotonicities. The
rearrangement delivers curves that improve upon the original estimates
and that satisfy the natural monotonicity requirement. We quantify
this improvement in the next subsection.









\begin{figure}
\includegraphics[width = \textwidth, height = \textwidth]{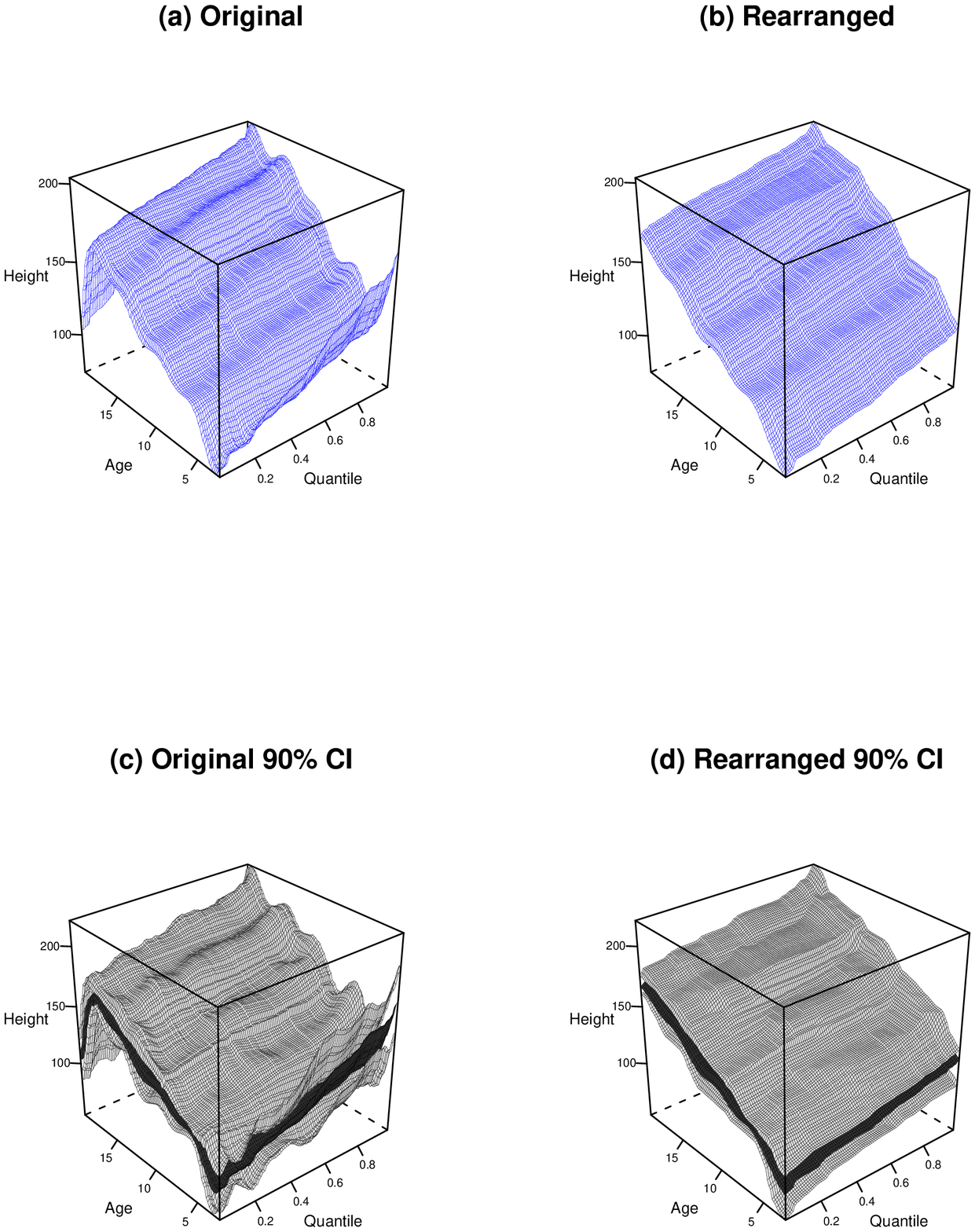}
\caption{Fourier series point and interval estimates of the
conditional quantile process of height given age and their
increasing rearrangements. Panels (a) and (b) plot original estimate
and its average multivariate rearrangement.  Panels (c) and (d) plot
original and rearranged 90\% confidence intervals. Original
confidence interval obtained by bootstrap with 200 repetitions.}
\label{Fig:7}
\end{figure}

Figure \ref{Fig:7} (a,b) illustrates the multivariate rearrangement
of the conditional quantile process  along both the age and the
quantile index arguments.  We plot, in three dimensions, the
original estimate and its average multivariate rearrangement (the
average of the age-quantile and quantile-age rearrangements). We
focus on the Fourier series estimates, which have the most severe
non-monotonicity problems. Analogous figures for the other
estimation methods are given in an MIT working paper containing an
extended version of this article. We see that the estimated quantile
process is non-monotone in age and in the quantile index at extremal
values of this index. The average multivariate rearrangement fixes
the non-monotonicity problem delivering an estimate of the quantile
process that is monotone in both the age and the quantile index.
Furthermore, by the theoretical results of the paper, the
multivariate rearranged estimates necessarily improve upon the
original estimates.

\begin{figure}
\includegraphics[width = \textwidth, height = \textwidth]{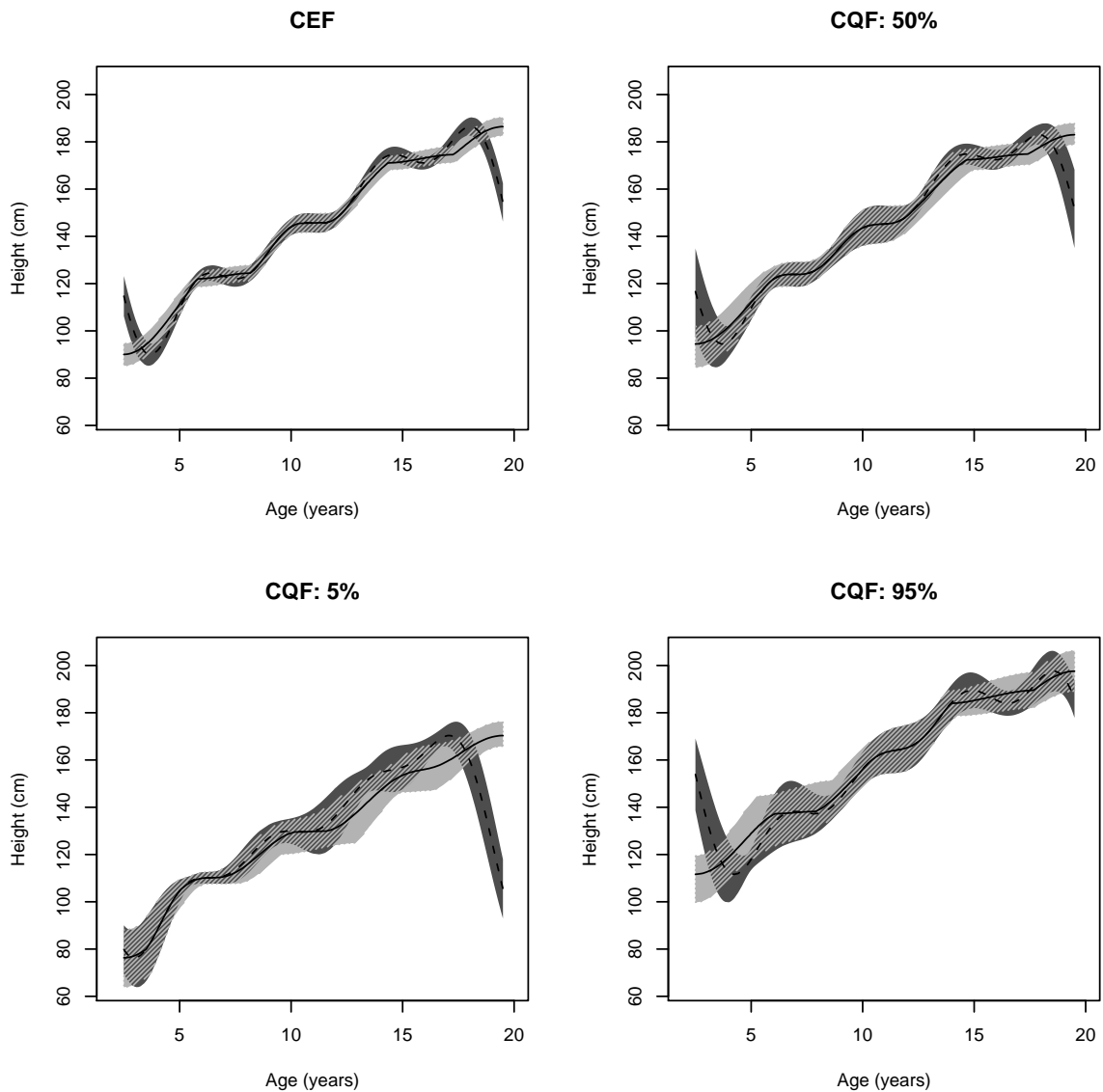}
\caption{\label{Fig:10} 90\% confidence intervals for conditional
expectation function (CEF), and 5\%, 50\% and 95\% conditional
quantile functions (CQF) of height given age and their increasing
rearrangements. Estimates are based on Fourier series and confidence
bands are obtained by bootstrap with 200 repetitions. Dark bands are
the original confidence intervals and light bands are the rearranged
confidence intervals.}
\end{figure}

In Figures \ref{Fig:10} and \ref{Fig:7} (c,d), we plot original and
rearranged 90\% simultaneous confidence intervals. Fig. \ref{Fig:10}
shows the intervals for the conditional expectation function and for
the conditional 5\%, 50\%, and 95\% quantile functions, based on
Fourier series estimates. We obtain the original intervals of the
form (\ref{CI2}) using the bootstrap with 200 repetitions
 to estimate the standard errors and critical values \cite{Hall1993}. We then obtain the
rearranged confidence intervals by rearranging the lower and upper
end-point functions of the initial confidence intervals, following
Section 3.  In Fig. \ref{Fig:7} (c,d), we plot the original and the
rearranged 90\% simultaneous confidence intervals for the entire
conditional quantile process, based on the Fourier series estimates.
The rearranged confidence intervals correct the non-monotonicity of
the original confidence intervals and reduce their integrated $L^p$
length.

\subsection{Monte-Carlo illustration}

In the following Monte Carlo experiment we quantify the improvement
in the point and interval estimation that rearrangement can provide
relative to the original estimates. We also compare it to
isotonization and to its convex combinations with isotonization. Our
experiment uses a model, described in detail in the Appendix, that
mimics the empirical application very closely.  This model implies a
true conditional expectation function and quantile process that are
monotone in age and in the quantile index.

In Table \ref{table1} we report the average $L^p$ errors, for
$p=1,2,$ and $\infty$, for the original estimates of the conditional
expectation function.  We also report the relative efficiency of the
rearranged estimates, measured as the ratio of the average error of
the rearranged estimate to the average error of the original
estimate; together with relative efficiencies for alternative
approaches based on isotonization of the original estimates
\cite{Mammen1991} and on averaging the rearranged and isotonized
estimates. For regression splines, we also consider the one-step
monotone regression splines \cite{Ramsay1998}.

For all of the methods and norms considered, the rearranged curves
estimate the target function more accurately than the original
curves. There is no uniform winner between rearrangement,
isotonization, and the average of the two, which is consistent with
the analysis of Section 2.4. For example, the rearrangement
outperforms the other methods for kernel, local linear and splines,
but performs worse than the average for Fourier in some norms. In
numerical results not reported, we find that rearrangement performs
worse than isotonization for global polynomials. This and other
methods are available in the MIT working paper. For regression
splines, the performance of the rearrangement is comparable to the
computationally more intensive one-step monotone splines procedure.

\begin{table}
\caption{$L^p$ Estimation Errors of Original, Rearranged,
Isotonized, Average Rearranged-Isotonized, and Monotone Estimates of
the Conditional Expectation Function, for $p=1,2$, $\infty$.}{
\begin{tabular}{c c c c c c c c c c c }
 $p$  & $\ \ L^p_O \ \ $& $L^p_R/L^p_O$ &$L^p_I/L^p_O$
&$L^p_{(R+I)/2}/L^p_O$ &$L^p_{M}/L^p_O$ & &$\ \ L^p_O \ \
$&$L^p_R/L^p_O$ &$L^p_I/L^p_O$ &$L^p_{(R+I)/2}/L^p_O$
\\
 \multicolumn{1}{l}{ }&\multicolumn{5}{c}{Kernel} &
\multicolumn{1}{l}{ }& \multicolumn{4}{c}
{ Local linear} \\
$1$      & 1$\cdot$00 & 0$\cdot$97 &0$\cdot$98 &0$\cdot$98& --  & & 0$\cdot$79& 0$\cdot$96& 0$\cdot$97& 0$\cdot$96\\
$2$      & 1$\cdot$30 & 0$\cdot$98 &0$\cdot$99 &0$\cdot$98& --  & & 0$\cdot$99& 0$\cdot$96& 0$\cdot$97& 0$\cdot$97\\
$\infty$ & 4$\cdot$54 & 0$\cdot$99 &1$\cdot$00 &1$\cdot$00& --  & & 2$\cdot$93& 0$\cdot$95& 0$\cdot$95& 0$\cdot$95\\
\multicolumn{1}{l}{ }&\multicolumn{5}{c}{Regression splines} &
\multicolumn{1}{l}{ }& \multicolumn{4}{c}
{Fourier} \\
$1$      & 0$\cdot$87 &0$\cdot$93 &0$\cdot$95 &0$\cdot$94 &0$\cdot$99 & & 6$\cdot$57 &0$\cdot$49 &0$\cdot$59 &0$\cdot$40\\
$2$      & 1$\cdot$09 &0$\cdot$93 &0$\cdot$95 &0$\cdot$94 &0$\cdot$99 & & 10$\cdot$8 &0$\cdot$35 &0$\cdot$45 &0$\cdot$30\\
$\infty$ & 3$\cdot$68 &0$\cdot$85 &0$\cdot$88 &0$\cdot$86 &0$\cdot$84 & & 48$\cdot$9 &0$\cdot$16 &0$\cdot$34 &0$\cdot$20\\
\end{tabular}}\label{table1}
\begin{flushleft}
\footnotesize{$L^p_O$, $L^p_{R}$,$L^p_{I}$,
 $L^p_{(R+I)/2}$, and $L^p_{M}$ are  average $L^p$ errors of the original, rearranged, isotonized, the average rearranged-isotonized, and monotone regression splines estimates; they are  computed as the Monte Carlo average of
$\{ \int_{\mathcal{X}} | \bar{f}(x) - f_{0}(x) |^{p} dx\}^{1/p},$
where $f_{0}$ is the target  and $\bar{f}$ an estimate.}
\end{flushleft}
\end{table}

In Table \ref{table2} we report the average $L^p$ errors for the
original estimates of the conditional quantile process. We also
report the ratio of the average error of the multivariate rearranged
estimate, with respect to the age and quantile index arguments, to
the average error of the original estimate; together with the same
ratios for isotonized and average rearranged-isotonized estimates.
We obtain the multivariate isotonized estimates by sequentially
applying the univariate isotonization to each argument, and then
averaging for the two possible orderings age-quantile and
quantile-age. For all the methods and norms considered, the
multivariate rearranged curves estimate the target function more
accurately than the original curves. There is again no uniform
winner between rearrangement, isotonization, and their average.

\begin{table}
\caption{$L^p$ Estimation Errors of Original, Rearranged,
Isotonized, and Average Rearranged-Isotonized Estimates of the
Conditional Quantile Process, for $p=1,2,$ and $\infty$. }{
\begin{tabular}{c c c c c c c c c c c }
  $p$ &  & $\ \ L^p_O \ \ $ &$L^p_{R}/L^p_O$
&$L^p_{I}/L^p_O$ &$L^p_{(R+I)/2}/L^p_O$ & &$\ \ L^p_O \ \ $ &
$L^p_{R}/L^p_O$ &$L^p_{I}/L^p_O$ &$L^p_{(R+I)/2}/L^p_O$
\\
 \multicolumn{2}{l}{ }&\multicolumn{4}{c}{Kernel} &
\multicolumn{1}{l}{ }& \multicolumn{4}{c}
{Local linear} \\
$1$     & & 1$\cdot$49  &0$\cdot$95 &0$\cdot$97 &0$\cdot$96 & & 1$\cdot$21& 0$\cdot$91& 0$\cdot$93& 0$\cdot$92\\
$2$     & & 1$\cdot$99  &0$\cdot$96 &0$\cdot$98 &0$\cdot$97 & & 1$\cdot$61& 0$\cdot$91& 0$\cdot$93& 0$\cdot$92\\
$\infty$& & 13$\cdot$7 & 0$\cdot$92 &0$\cdot$97 &0$\cdot$94 & & 12$\cdot$3& 0$\cdot$84& 0$\cdot$87& 0$\cdot$85\\
\multicolumn{2}{l}{ }&\multicolumn{4}{c}{Regression splines} &
\multicolumn{1}{l}{ }& \multicolumn{4}{c}
{Fourier} \\
$1$     & & 1$\cdot$33  &0$\cdot$90 &0$\cdot$93 &0$\cdot$91 & & 6$\cdot$72  &0$\cdot$62 &0$\cdot$77 &0$\cdot$64\\
$2$     & & 1$\cdot$78  &0$\cdot$90 &0$\cdot$92 &0$\cdot$90 & & 13$\cdot$7  &0$\cdot$39 &0$\cdot$58 &0$\cdot$44\\
$\infty$& & 16$\cdot$9 & 0$\cdot$72 &0$\cdot$76 &0$\cdot$73 & &
84$\cdot$9 & 0$\cdot$26 &0$\cdot$47
&0$\cdot$36\\
\end{tabular}}\label{table2}
\begin{flushleft}
\footnotesize{$L^p_O$, $L^p_{R}$,$L^p_{I}$,
 and $L^p_{(R+I)/2}$ are the average $L^p$ errors of the original, multivariate rearranged,  multivariate isotonized, the multivariate average rearranged-isotonized estimates; they are computed as the Monte Carlo averages of
$\{\int_{\mathcal{U}} \int_{\mathcal{X}} | \bar{f}(u, x) - f_{0}(u,
x) |^{p}dx du\}^{1/p},$ where $f_{0}$ is the target and $\bar{f}$ an
estimate.}
\end{flushleft}

\end{table}

Table \ref{table3} reports Monte Carlo coverage frequencies and
integrated lengths for the original and monotonized 90\% confidence
bands for the conditional expectation function.  For a measure of
length, we used the integrated $L^p$ length, as defined in
Proposition 3, with  $p = 1,2,$ and $\infty$. We construct the
original confidence intervals of the form specified in equations
(\ref{CI2}) by obtaining the pointwise standard errors of the
original estimates using the bootstrap with 200 repetitions, and
calibrate the critical value so that the original confidence bands
cover the entire true function with the exact frequency of 90\%. We
 construct monotonized confidence intervals by applying
rearrangement, isotonization,  and a rearrangement-isotonization
average to the end-point functions of the original confidence
intervals, as proposed in Section 3.   In all cases the
rearrangement and other monotonization methods increase the coverage
of the confidence intervals while reducing their length.  In
particular, we see that monotonization increases coverage especially
for the local estimation methods, whereas it reduces length most
noticeably for the global estimation methods. For the most
problematic Fourier estimates, there are large increases in coverage
and reductions in length.

\begin{table}
\caption{Coverage (\%) and Integrated Lengths of Original,
Rearranged, Isotonized, and Average Rearranged-Isotonized 90\%
Confidence Intervals for the Conditional Expectation Function.}{
\begin{tabular}{c c c c c c c c c c c}
  \multicolumn{1}{l}{Interval
}&\multicolumn{1}{c}{Cover} & \multicolumn{4}{c}{Length} &
\multicolumn{1}{c}{Cover} &
\multicolumn{4}{c}{Length}  \\
 & & $\ \ L^1 \ \ $ &$L^1/L^1_O$
&$L^2/L^2_O$ &$L^{\infty}/L^{\infty}_O$ & &$\ \ L^1 \ \
$&$L^1/L^1_O$ &$L^2/L^2_O$ &$L^{\infty}/L^{\infty}_O$
\\
 \multicolumn{1}{l}{ }&\multicolumn{5}{c}{Kernel} &
 \multicolumn{5}{c}
{Local linear} \\
$O$&        90 & 8$\cdot$80 &   &   &       & 90 & 8$\cdot$63 &   &   & \\
$R$&        96 & 8$\cdot$79 & 1 & 1 & 0$\cdot$99  & 96 & 8$\cdot$63 & 1 & 1 & 0$\cdot$97\\
$I$&        94 & 8$\cdot$80 & 1 & 1 & 0$\cdot$99  & 94 & 8$\cdot$63 & 1 & 1 & 0$\cdot$98\\
$(R+I)/2$&  95 & 8$\cdot$80 & 1 & 1 & 0$\cdot$99  & 95 & 8$\cdot$63 & 1 & 1 & 0$\cdot$97\\
\multicolumn{1}{l}{ }&\multicolumn{5}{c}{Regression splines} &
 \multicolumn{5}{c}
{Fourier} \\
$O$&        90 & 6$\cdot$32 &   &   &       & 90  & 24$\cdot$91 &      &      &     \\
$R$&        91 & 6$\cdot$32 & 1 & 1 & 1     & 100 & 24$\cdot$52 & 0$\cdot$98 & 0$\cdot$94 & 0$\cdot$63\\
$I$&        91 & 6$\cdot$32 & 1 & 1 & 1     & 100 & 24$\cdot$91 & 1    & 0$\cdot$97 & 0$\cdot$69\\
$(R+I)/2$&  91 & 6$\cdot$32 & 1 & 1 & 1     & 100 & 24$\cdot$71 & 0$\cdot$99 & 0$\cdot$95 & 0$\cdot$65\\
\end{tabular}}\label{table3}
\begin{flushleft}
\footnotesize{$O$, $R$, $I$, and $(R+I)/2$ refer to original,
rearranged, isotonized, and average rearranged-isotonized confidence
intervals. Coverage probabilities (Cover) are for the entire
function.}
\end{flushleft}

\end{table}

\appendix
\vspace{-.1in}

\section*{Acknowledgment}
We would like to thank the editor, the associate editor, many
referees of the journal, M. Cohen, H. Dette, E. Gallagher, W.
Graybill, P. Groneboem, R. Guiteras, X. He, R. Koenker, C. Manski,
I. Molchanov, W. Newey, S. Portnoy, A. Simsek, J. Wellner and
participants of many seminars and conferences for comments that
helped to considerably improve the paper. Chernozhukov and
Fern\'andez-Val gratefully acknowledge research support from the
NSF. Galichon's research is partly supported by chaire
X-Dauphine-EDF-Calyon ``Finance et D\'eveloppement Durable''.

\vspace{-.1in}

\appendix

\section{Proofs of Propositions}

\subsection*{Proof of Proposition 1}

Proof of Part 1. This follows in part the strategy in Lorentz's
(1953) proof. We assume first  that the functions $\hat f$ and $f_0$
are step functions, constant on intervals $((s-1)/r, s/r]$,
$s=1,\ldots,r$. For each step function $f$ with $r$ steps we
associate an $r$-vector $f$ whose $s$-th element, denoted $f_{s}$,
equals to the value of function $f$ on the $s$-th interval, and vice
versa. Let us define the sorting operator $S$ acting on vectors (and
functions) $f$ as follows. Let $k$ be an integer in $1,\ldots,r$
such that $f_k > f_m$ for some $m>k$. If $k$ does not exist, set $Sf
= f$. If $k$ exists, set $Sf$ to be a $r$-vector with the $k$-th
element equal to $f_m$, the $m$-th element equal to $f_k$, and all
other elements equal to the corresponding elements of $f$. Finally,
given a vector $Sf$ there is a step function $Sf$ associated to it,
as stated above.

For any submodular function $L: \Bbb{R}^2 \to \Bbb{R}_+$, by $f_k
\geq f_m$, $f_{0m} \geq f_{0k}$ and the definition of the
submodularity, $L(f_m, f_{0k}) + L (f_k, f_{0m}) \leq L(f_k, f_{0k})
+ L(f_m, f_{0m}).$  A simple geometric illustration for this
property is given in Figure \ref{Fig:prop1}.   Therefore, conclude
that $ \int_{\X} L\{S\hat f(x), f_0(x)\} dx \leq \int_{\X} L\{ \hat
f(x), f_0(x)\} dx,$ using that  we integrate step functions.
Applying the sorting operator a sufficient finite number of times to
$\hat f$, we obtain a completely sorted, that is, rearranged, vector
$\hat f^*$. Thus, we can express $\hat f^*$ as  $\hat f^*= S \ldots
S \hat f$, where the operator $S$ is applied finitely many times. By
repeating the argument above, each application weakly reduces the
estimation error. Therefore,
\begin{equation}\label{maininequality}
\int_\X L\{\hat f^*(x), f_0(x)\} d x \leq \int_{\X} L\{ S \ldots
S\hat f(x) , f_0(x) \} dx \leq \int_{\X} L\{\hat f(x), f_0(x)\} dx.
 \end{equation}

Next we extend this result to general measurable functions $\hat f$
and $f_0$ mapping $[0,1]$ to $K$, where $f_0$ is a quantile
function. Take a subsequence of  bounded step functions $\hat
f^{(q)}$ and $f_0^{(q)}$, with $f_0^{(q)}$ being quantile functions,
converging to $\hat f$ and $f_0$ almost everywhere as index $q \to
\infty$ along an increasing sequence of integers. The almost
everywhere convergence of $\hat f^{(q)}$ to $\hat f$ implies the
almost everywhere convergence of its quantile function $\hat
f^{*(q)}$ to the quantile function of the limit, $\hat f^*$
(\citeasnoun{vaart:text}, p. 305). Since (\ref{maininequality})
holds for each $q$ along the subsequence, the dominated convergence
theorem implies that (\ref{maininequality}) also holds for the
general case.

It remains to show the existence of the subsequence in the preceding
paragraph. Using series expansion in the Haar basis, any function in
$L^2[0,1]$ can be approximated in $L^2$ norm  by a sequence of
$r$-step functions, where $r=2^j$ and $j =1,\ldots,\infty$
(\citeasnoun{pollard:measure}, p. 305) . Hence there is a
subsequence of step functions $ \hat f^{(r)}$ and $f^{(r)}_0$
converging to $\hat f$ and $f_0$ in $L^2$ norm; the functions in the
subsequence necessarily take values in $K$; by
\citeasnoun{pollard:measure}, p. 38, we can extract a further
subsequence $\hat f^{(q)}$ and $f^{(q)}_0$, with $q$ running over an
increasing sequence of integers, converging to $\hat f$ and $f_0$
almost everywhere. Finally, replace $f^{(q)}_0$ by their quantile
functions, i.e., rearrangements, which retain the almost everywhere
convergence property to $f_0$ by \citeasnoun{vaart:text}, p. 305.
\qed

Proof of Part 2. Consider the step functions, as defined in the
proof of Part 1.  By setting $r$ sufficiently large, we  can take
them to satisfy the following hypotheses: there exist regions $\X_0$
and $ \X'_0$, each of measure greater than $\delta>0$, such that for
all $x \in \X_0$ and $x' \in \X_0'$, we have that (i) $x'
> x$, (ii) $\hat f(x) > \hat f(x') + \epsilon$, and (iii)
$f_0(x') > f_0(x) + \epsilon$, for $\epsilon>0$ specified in the
proposition.   For any strictly submodular function $L: \Bbb{R}^2
\to \Bbb{R}_+$ we have that $\eta = \inf \{ L(v',t) + L(v,t') -
L(v,t) - L(v',t') \} >0, $ where the infimum is taken over all $v,
v', t, t'$ in the set $K$ such that $v' \geq v + \epsilon$ and $t'
\geq t + \epsilon$.  We can begin sorting by exchanging an element $\hat f(x)$, $x \in
\X_0$, of $r$-vector $\hat f$ with an element $\hat f(x')$, $x' \in
\X_0'$, of $r$-vector
 $\hat f$. This induces a sorting gain of at
least $\eta$ times $1/r$.  The total mass of points that can be
sorted in this way is at least $\delta$. We then proceed to sort all
of these points in this way, and then continue with the sorting of
other points. After the sorting is completed, the total gain from
sorting is at least $\delta \eta$. That is, $ \int_{\X} L\{\hat
f^*(x), f_0(x)\}dx  \leq \int_{\X} L\{\hat f(x), f_0(x)\} dx -
\delta \eta.$

We then extend this inequality to the general measurable functions
exactly as in the proof of Part 1. \qed

\subsection*{Proof of Proposition 2}

Proof of Part 1.  We prove the claim by induction. It is true
for $d=1$ by $\hat f^*$ being a quantile function. Suppose the claim is true in $d-1 \geq 1$
dimensions. If so, then $x_{-j}\mapsto \bar f(x_j, x_{-j})$, obtained
from the original estimate $\hat f$ after applying the
rearrangement to all arguments $x_{-j}$ of $x$, except for the
argument $x_j$, must be weakly increasing in $x_{-j}$ for each
$x_j$. Thus, for any $x_{-j}' \geq x_{-j}$ and  $X_j \sim
U[0,1]$, we have
\begin{equation}\label{dom1}
\bar f(X_j, x_{-j}') \geq \bar f(X_j, x_{-j}).
 \end{equation}
Therefore, the random variable on the left of (\ref{dom1}) dominates
the random variable on the right of (\ref{dom1}) in the stochastic
sense. Therefore, the quantile function of the random variable on
the left dominates the quantile function of the random variable on
the right, namely $
\bar f^*_j(x_j, x_{-j}') \geq \bar f^*_j(x_j, x_{-j}) \text{ for
each }  x_j \in \X = [0,1].$
Moreover, for each $x_{-j}$, the function $x_j \mapsto \bar
f^*_j(x_j, x_{-j})$ is weakly increasing by virtue of being a
quantile function.  We conclude therefore that $x\mapsto \bar
f_j^*(x)$ is weakly increasing in all of its arguments at all points
$x \in \X^d$.  The claim of Part 1 of the Proposition now follows by
induction. \qed

Proof of Part 2 (a).  By Proposition 1, we have that for each
$x_{-j}$,
 \bsnumber\label{weak} \int_{\X}
\left | \bar f_j^*(x_j, x_{-j}) - f_0(x_j, x_{-j}) \right |^p d x_j
\leq \int_{\X} \left | \bar f (x_j, x_{-j}) - f_0(x_j, x_{-j})
\right |^p d x_j.
 \end{split}\end{align}
Now, the claim follows by integrating with respect to $x_{-j}$ and
taking the $p$-th root of both sides. For $p = \infty$, the claim
follows by taking the limit as $p \to \infty$. \qed

Proof of Part 2 (b).  We first apply the inequality of Part 2(a) to
$\bar f(x) = \hat f(x)$, then to $\bar f(x) = R_{\pi_d} \hat f(x),$
then to $\bar f(x) = R_{\pi_{d-1}} R_{\pi_d} \hat f(x)$, and
so on.   In doing so, we recursively generate a sequence of weak
inequalities that imply the inequality (\ref{P22b}) stated in the
Proposition. \qed

Proof of Part 3 (a).  For each $x_{-j} \in \X^{d-1}\setminus
\X_{-j}$, by Part 2(a), we have the weak inequality (\ref{weak}),
and for each $x_{-j} \in \X_{-j}$, by the inequality for the
univariate case stated in Proposition 1 Part 2, we have the strong
inequality \bsnumber\label{strong} \int_{\X} \left | \bar f_j^*(x_j,
x_{-j}) - f_0(x_j, x_{-j}) \right |^p d x_j \leq \int_{\X} \left |
\bar f (x_j, x_{-j}) - f_0(x_j, x_{-j})  \right |^p d x_j - \eta_p
\delta,
 \end{split}\end{align}
where $\eta_p$ is defined in the same way as in Proposition 1.
Integrating the weak inequality (\ref{weak}) over $x_{-j} \in
\X^{d-1}\setminus \X_{-j}$, of measure $1- \nu$, and the strong
inequality (\ref{strong}) over $\X_{-j}$, of measure $\nu$, we
obtain \bsnumber\label{strong2} \int_{\X^d} \left | \bar f_j^*(x) -
f_0(x) \right |^p d x \leq \int_{\X^d} \left | \bar f (x) - f_0(x)
\right |^p d x - \eta_p \delta \nu.
 \end{split}\end{align}
 The claim now follows. \qed

 Proof of Part 3 (b).  As in Part 2(a), we can
 recursively obtain a sequence of weak inequalities
 describing the improvements in estimation error
 from rearranging sequentially with respect to the individual
 arguments.  Moreover, at least one of the inequalities
 can be strengthened to be of the form stated in
 (\ref{strong2}), from the assumption of the claim.
 The resulting system of
 inequalities yields the inequality (\ref{P23b}),
 stated in the proposition. \qed

 Proof of Part 4.  This part follows from homogeneity and subadditivity of the $L_p$ norm.  \qed

\subsection*{Proof of Proposition 3}

Proof of Part 1.  The monotonicity follows from Proposition 2.
 The rest of the proof relies on establishing the order-preserving property
of the $\pi$-rearrangement operator:  for any measurable functions
$g, m : \mathcal{X}^d \to \Bbb{R}$, we have that
$g(x) \leq m(x)$ for all  $x \in \mathcal{X}^d$
implies  $g^*(x) \leq  m^*(x)$ for all $x \in \mathcal{X}^d.$
Given the property we have that $\ell(x) \leq f(x) \leq u(x)  \text{ for all } x \in \mathcal{X}^d
\text{ implies } \ell^*(x) \leq  f^*(x)  \leq u^*(x) \text{ for all
} x \in \mathcal{X}^d,$ which verifies the claim of the first part.  The claim also extends
to the average multivariate rearrangement, since averaging preserves
the order-preserving property.

It remains to establish the order-preserving property for
$\pi-$rearrangement, which we do by induction. We first note that in
the univariate case, when $d=1$, order preservation is obvious from
the rearrangement being a quantile function: the random variable
$m(X)$, where $X \sim U[0,1]$, dominates the random variable $g(X)$
in the stochastic sense, hence the quantile function $m^*(x)$  of
$m(X)$ must be weakly greater than the quantile function $g^*(x)$ of
$g(X)$ for each $x \in \mathcal{X}$. We then extend this to the
multivariate case by induction: Suppose the order-preserving
property is true for any $d-1 \geq 1$.  If so, then, for each $x_j
\in \mathcal{X}$ and $x_{-j} \in \mathcal{X}^{d-1}$, $ g(x_j,
x_{-j}) \leq m(x_j, x_{-j}) \text{ implies }  \bar g (x_j, x_{-j})
\leq  \bar m (x_j, x_{-j}),$ where $\bar g$ and $\bar m$ are
multivariate rearrangements of $x_{-j} \mapsto  g (x_j, x_{-j})$ and
$x_{-j} \mapsto m (x_j, x_{-j})$ with respect to $x_{-j}$, holding
$x_j$ fixed. Now apply the order-preserving property of the
univariate rearrangement to the univariate functions $x_{j} \mapsto
\bar g (x_j, x_{-j})$ and $x_{j} \mapsto \bar m (x_j, x_{-j})$,
holding $x_{-j}$ fixed, for each $x_{-j}$, to conclude that the
order-preserving property holds for dimension $d$.\qed

Proof of Part 2. As stated in the text, the weak inequality follows
from Lorentz (1953).  For completeness we only briefly note that the
proof follows similarly to the proof of Proposition 1. Indeed, we
can start with step functions $\ell$ and $u$ and work with their
equivalent vector representations $\ell$ and $u$. Then we apply the
sorting operator $S_2$ to the pair of $r$-vectors $(\ell, u)$
defined as $S_2(\ell, u)=(S\ell, S'u)$, where $S$ is the sorting
operator on the vector $\ell$ defined in the proof of Proposition 1,
and $S'$ is a subordinated sorting operator on the vector $u$
defined by two conditions: (1) if $S$ exchanges the $k$-th and
$m$-th elements of $\ell$, where $m>k$, then,  if $u_k > u_m$, $S'$
also exchanges the $k$-th and $m$-th elements of $u$ , and if $u_k
\leq u_m$, $S'$ leaves all elements of $u$ unchanged; and (2) if $S$
exchanges no elements of $\ell$, i.e., $S\ell = \ell$, then $S'$ is
simply the unrestricted $S$ operator as defined in the proof of
Proposition 1, i.e., $S'=S$. By the definition of submodularity
(\ref{submodular}), each application of $S_2$ weakly reduces
submodular discrepancies between vectors, so that the pairs of
vectors in the sequence $ \{ (\ell, u), S_2(\ell, u),\ldots,
S_2\ldots S_2(\ell, u),  (\ell^*, u^*) \}$ become progressively
weakly closer to each other, and the sequence can be taken to be
finite, where the last pair is the rearrangement $(\ell^*, u^*)$ of
vectors $(\ell,u)$.  The inequality extends  to general bounded
measurable functions by passing to the limit using a similar
argument to the proof of Proposition 1. The extension of the proof
to the multivariate case follows by induction on the dimension, as
in the proof of Proposition 2.  \qed

Proof of Part 3. Finally, the proof of strict inequality in the
univariate case is similar to the proof of Proposition 2,
  using the fact that for strictly submodular functions $L: \Bbb{R}^2 \mapsto \Bbb{R}_+$ we have that $\eta = \inf \{ L(v',t) + L(v,t') -
L(v,t) - L(v',t') \} >0, $ where the infimum is taken over all $v,
v', t, t'$ in the set $K$ such that $v' \geq v + \epsilon$ and $t'
\geq t + \epsilon$ or  such that $v \geq v' + \epsilon$ and $t \geq
t' + \epsilon$.  The extension of the strict inequality to the
multivariate case follows exactly as in the proof of Proposition 2.
\qed

\section{Design of the Monte-Carlo experiment}

The outcome variable $Y$ equals a location function plus a
disturbance $\epsilon$, $Y = Z(X)'\beta + \epsilon,$ and the
disturbance is independent of the regressor $X$. The vector $Z(X)$
includes a constant and a piecewise linear transformation of the
regressor $X$ with three changes of slope, namely $Z(X) =(1, X, 1\{
X > 5 \} (X - 5), 1\{ X > 10 \} (X - 10), 1\{ X > 15 \}(X - 15) )$.
This design implies the conditional expectation function $E[Y\mid X]
= Z(X)'\beta$, and the conditional quantile function  $Q_{Y}[u\mid
X] = Z(X)'\beta + Q_\epsilon(u).$ We select the parameters of the
design to match the growth charts example. Thus, we set the
parameter $\beta$ equal to the ordinary least squares estimate
obtained in the growth chart data, namely ($71.25$, $8.13$, $-2.72$,
$1.78$, $-6.43$). This parameter value and the location
specification imply a model for the conditional expectation function
and quantile process that is monotone for ages 2-20. To generate the
values of the dependent variable, we draw disturbances from a normal
distribution whose mean and variance match those of the estimated
residuals, $\epsilon = Y - Z(X)'\beta$. We fix the values of the
regressor $X$ to be the observed values of age in the data. In each
replication, we estimate the target functions  using the
nonparametric methods described in Section 4.1. The total number of
replications is 1000. All computations were carried out using the
software R \cite{R2007}, the quantile regression package quantreg,
and the functional data analysis package fda. The rearrangement
method developed in this paper is available in the package
rearrangement for R.

\vspace{-.2in}

\bibliographystyle{econometrica}
\bibliography{c:/aaa/biblio/my}

\end{document}